\documentclass[a4paper,12pt]{article}
\usepackage[english]{babel}
\usepackage[T2A]{fontenc}
\usepackage{amsthm}
\usepackage[tbtags]{amsmath}
\usepackage{amsfonts,amssymb}
\begin{document}

\begin{center}
\textbf{Generalized Bassian Modules\\  over Non-primitive Dedekind Prime Rings}

\textbf{Askar Tuganbaev}
\end{center}
Department of Higher Mathematics, National Research University MPEI,
Krasnokazarmennaya 14, Moscow, 111250, Russia; tuganbaev@gmail.com

\textbf{Abstract.} A right $A$-module $M$ is said to be generalized bassian if the existence of an injective homomorphism $M\to M/N$ for some submodule $N$ of $M$ implies that $N$ is a direct summand of $M$. We describe singular  generalized bassian modules over non-primitive Dedekind prime rings.\\
The study is supported by grant of Russian Science Foundation.

\textbf{Keywords:}  bassian module, generalized bassian module

Mathematics Subject Classification 2020: 16D10, 16D40, 16D60, 20K10, 20K20, 20K21

\section{Introduction}\label{sec1}
We only consider associative unital non-zero rings and unitary modules. The phrases of the kind <<a noetherian ring $A$>> mean that both modules $A_A$ and $_AA$ are noetherian. The word <<$A$-module>> usually means <<right $A$-module>>. 
The notions that are not defined in the paper are standard; e.g., see \cite{Fai76}-- \cite{Kas82}, \cite{Lam99}.

A module $M$ is said to be \textbf{generalized bassian} if the existence of a monomorphism $M\to M/N$, for some submodule $N$ of $M$, forces $N$ to be a direct summand of $M$. \\
A module $M$ is said to be \textbf{bassian} if $M$ is not isomorphic to a submodule of any its proper homomorphic image; equivalently, $N= 0$ for any submodule $N$ of $M$ such that there exists a monomorphism $M\to M/N$. Clearly, every bassian module is generalized bassian. 

\textbf{Remark 1.1.} If $A$ is a ring and $M$ is a direct sum of infinitely many isomorphic simple $A$-modules, then it is clear that $M$  is a generalized bassian module which is not bassian. In particular,  any elementary non-noetherian primary abelian group is a generalized bassian group which is not bassian. Theorem 1 of \cite{Tug26} describes all singular\footnote{It is well known that torsion abelian groups coincide with singular $\mathbb{Z}$-modules.} Bassian modules over non-primitive Dedekind prime rings (examples of such rings include all rings of algebraic integers and matrix rings over them; in particular,  the ring $\mathbb{Z}$ is a non-primitive Dedekind prime ring).  Some definitions used to formulate Theorem 1.2 below are given in the Introduction, after Theorem 1.2. Some other necessary definitions are given in Section 2.

 All bassian abelian groups are described in \cite[Main Theorem]{CheDG21}. The description of  generalized bassian abelian groups follows from  \cite[Theorem 2.10, Corollary 3.14]{CheDG22},  \cite[Corollary 3.6]{DanK25}, and  \cite[Corollary 4.4]{DasK26}. 

The main result of  the paper is the following Theorem 1.2.

\textbf{Theorem 1.2. } If $A$ is a non-primitive Dedekind prime ring and $M$ is a singular right $A$-module with primary components $M_i$, $i\in I$, then the following conditions are equivalent.
 \begin{enumerate} 
 \item[ 1)]  
 $M$ is a generalized Bassian module. 
 \item[2)] 
 $M$ is a direct sum of a Bassian module and a semisimple module. 
 \item[3)] 
Every primary component $M_i$ of $M$ is a direct sum of a noetherian module and a semisimple module. \end{enumerate}

A module $M_A$ is said to be \textbf{projective} if it satisfies the following equivalent conditions: \textbf{1)} for any module $X_A$, every epimorphism $h\colon X\to \overline{X}$, and an arbitrary homomorphism $\overline{f}\colon M\to \overline{X}$, there exists a homomorphism $f\colon M\to X$ with $\overline{f}=hf$; \textbf{2)} the module $M$ is isomorphic to a direct summand of some direct sum of copies of the module $A_A$. A module is said to be \textbf{hereditary} if all its submodules are projective. A module $M$ is said to be \textbf{noetherian} if all its submodules are finitely generated (equivalently, $M$ is a module with maximum condition on ascending chains of submodules). 

A ring is said to be \textbf{semiprime} (resp., \textbf{prime}) if for any its non-zero ideal $B$ the ideal $B^n$ (resp., the ideal $BC$ for any its non-zero ideal $C$) is non-zero. 
Let $A$ be a  ring which has a two-sided classical right and left ring of fractions $Q$ isomorphic to a matrix ring over a division ring. An ideal $B$ of $A$ is said to be \textbf{invertible} if there exists a sub-bimodule $B^{-1}$ of the bimodule ${}_AQ_A$ such that $BB^{-1}=B^{-1}B=A$. The maximal elements of the set of all proper invertible ideals of the ring $A$ are called \textbf{maximal invertible} ideals.
The set of all maximal invertible ideals of the ring $A$ is denoted by ${\mathcal P}(A)$. If $M$ is an $A$-module and $P\in {\mathcal P}(A)$, then the submodule $\{m\in M\;\vert \;mP^n=0, ~n=1,2, \ldots \}$ is called the \textbf{$P$-primary component} of the module $M$; it is denoted by $M(P)$. If $M=M(P)$ for some $P\in {\mathcal P}(A)$, then the module $M$ is said to be \textbf{primary} or \textbf{$P$-primary}.

A ring $A$ is called a \textbf{Dedekind prime} ring if $A$ is a noetherian prime ring\footnote{It is well known that  any noetherian prime ring has a two-sided classical right and left ring of fractions $Q$ isomorphic to a matrix ring over a division ring.} with classical ring of fractions $Q$ and every non-zero ideal of the ring $A$ is invertible in $Q$. A ring $A$ is said to be \textbf{non-primitive} if $A$ is not right or left primitive, i.e. $A$ has no faithful simple right or left modules. 
For a  module $M$, a submodule $X$ of $M$ is said to be \textit{essential} if $X\cap Y\ne 0$ for every non-zero submodule $Y$ of $M$. In this case, one says that $M$ is an essential extension of $X$. A submodule $X$ of $M$ is said to be \textbf{closed} in $M$ if $X$ has no proper essential extensions in $M$.  A module $M$ is said to be \textit{uniform} if the intersection of any pair of its non-zero submodules is non-zero, i.e. any non-zero submodule of $M$ is essential.  A module $M$ is said to be \textbf{uniserial} if the lattice of all its submodules is a chain, i.e. any two submodule of $M$ are comparable under inclusion (equivalently, any two cyclic submodules are comparable under inclusion). Every uniserial module is uniform and every uniform module is indecomposable. 

A ring $A$ is said to be \textbf{right} (resp., \textbf{left}) \textbf{bounded} if any its essential right (resp., \textbf{left}) ideal contains a non-zero ideal of the ring $A$. 

\textbf{Remark 1.3.} All Dedekind prime rings are hereditary noetherian prime (\textit{HNP}) rings. In \cite{Len73}, it is proved that any HNP ring $A$ is (right and left) primitive or (right and left) bounded and if $A$ is a primitive bounded ring, then $A$ is a simple artinian ring. Therefore, non-primitive HNP rings coincide with non-artinian bounded HNP rings.

For a subset $X$ of  a right (resp., left) $A$-module $M$, we  denote by $r(X)$  (resp., $\ell(X)$) the right (resp.,  left) annihilator of a $X$ in $A$, i.e. $r(X)=\{a\in A\,|\,Xa=0\}$ (resp., $\ell(X)=\{a\in A\,|\,aN=0\}$).  For a  right (resp., left) $A$-module $M$, we denote by $\text{Sing }M$ the set of all elements $m$ of $M$ such that $r(m)$ (resp., $\ell(m)$) is an essential right (resp., left) ideal of $A$. It is well known that $\text{Sing }M$ is a fully invariant submodule of $M$; it is called the \textbf{singular submodule} of $M$. If $\text{Sing }M=M$ (resp., $\text{Sing }M=0$), then the module $M$ is said to be \textbf{singular} (resp., \textbf{non-singular}). For a  module $M$, 2e denote by $Z_2(M)$ the \textbf{Goldie radical} of $M$, i.e. $Z_2(M)$ is the submodule of $M$ such that $Z_2(M)$ contains $\text{Sing }M$ and $Z_2(M)/\text{Sing }M)=\text{Sing }(M/\text{Sing }M$. A module $M$ is called a \textbf{Goldie-radical} module if $M=Z_2(M)$. 

If a module $M$ is an essential extension of  $\oplus_{j\in J}X_j$, where all $X_j$ are non-zero uniform modules and $\tau$ is the cardinality of the set $J$, then $\tau$ is called the \textbf{uniform dimension} or the \textbf{Goldie dimension} of $M$; it is denoted by $\text{u.dim }M$ or $\text{Gdim }M$. It is well known that the cardinal number $\text{u.dim }M$ is uniquely defined. A module $M$ is said to be \textbf{finite-dimensional} if $M$ has no submodules which are infinite direct sums of non-zero modules; in this case, its uniform dimension is finite.  A right finite-dimensional ring $A$ with the ascending chain condition on right annihilator ideals is called a \textbf{right Goldie} ring.

\textbf{Remark 1.4.} Let $A$ be a non-primitive HNP ring (e.g. we may take $A$ to be a non-primitive Dedekind prime ring), $M$ be a singular $A$-module and let $\{M_i\}_{i\in I}$ be the set of all primary components of $M$. Many properties of the singular module $M$ and its primary components are well known; for example, see \cite{Sin74}, \cite{Sin75}, \cite{Sin76}, \cite{Sin78}, \cite{Sin79}, \cite{SinT79}. These properties are similar to properties of primary components of torsion abelian groups. For example, $M=\oplus_{i\in I}M_i$,~ $X=\oplus_{i\in I}X_i$ for any submodule $X$ of the module $M$, where $X_i=X\cap M_i$ are primary components of the module $X$, $M/X=\oplus_{i\in I}(M_i/X_i)$, etc. Let $i,j\in I$, $i\ne j$, $X_i$ be a subfactor of $M_i$, and let $Y_j$ be a subfactor of $M_j$. Then there are no non-zero homomorphisms between $X_i$ and $Y_j$ (in particular, all primary components $M_i$ are fully invariant in $M$).

The symbols $\mathbb{Z}$ and $\mathbb{Q}$ denote the sets of integers and rational numbers, respectively.
\section{General Results}\label{sec2}

For a ring $A$, an element $a$ of $A$ is said to be \textbf{regular} if $r(a)=\ell(a)=0$. For a module $M_A$, we denote by $T(M)$ the set of all elements $m$ of $M$ such that $r(m)$ contains a regular element of $A$; this set is called the \textbf{torsion part} of $M$. A module $M$ is called a \textbf{torsion module} (resp., a \textbf{torsion-free module}) if $T(M)=M$ (resp., $T(M)=0$). 

\textbf{Remark 2.1 \cite[Example 2.1]{DasK26}, \cite[Remark 4]{Tug26}.} All noetherian modules are bassian. 

\textbf{Proposition 2.2  \cite[Lemma 2.1]{DasK26}, \cite[Remark 5]{Tug26}.}   If $M$ is a generalized bassian (resp., bassian) module and $M = X\oplus Y$, then $X$ is  a generalized bassian (resp., bassian) module.

\textbf{Proof.} Let $M$ be a generalized bassian module, $N$ be a submodule of $X$ and let $f\colon X\to X/N$ be a monomorphism. Then the monomorphism $f$ may be extended to a monomorphism $f\oplus 1_Y\colon X\oplus Y\to (X/N)\oplus Y=M/N$. Since $M$ is generalized bassian, there exists a submodule $N'$ of $M$ such that $M=N\oplus N'$, whence $X = (N\oplus N')\cap X = N\oplus (N'\cap X)$ by modularity. Therefore, $M$ is generalized bassian.

Let $M$ be a bassian module. Assume that the module $X$ is not Bassian. There exists a monomorphism $f\colon X \to X/X_1$ such that $X_1$ is a non-zero submodule of $X$. Then there exists a monomorphism $g = f \oplus 1_Y\colon M \to M/X_1= (X/X_1 )\oplus B$, where $X_1\ne 0$. This is a contradiction.\hfill$\square$

\textbf{Proposition 2.3 \cite[Propositions 3.1, 4.1, 4.4 and Corollary 4.2]{DasK26}.}\\
Let $M$ be a  generalized bassian module. 
\begin{enumerate}
\item[{\bf 1.}] 
If every homogeneous component of $\text{Soc }(M)$ is finitely generated, then the  module  $M$ is bassian. In particular, if $M$ has  a zero socle, then the  module  $M$ is bassian.
\item[{\bf 2.}] 
If  there is no a monomorphism from $M$ into any proper direct summand of $M$,  then $M$ is a bassian module. In particular, if either $M$ is indecomposable or any injective  endomorphism of $M$ is an automorphism  of $M$, then the  module $M$ is bassian.
\end{enumerate}

\textbf{Proposition 2.4.}  
Let $M$ be a module and let $X$ be a submodule of $M$. 
\begin{enumerate}
\item[{\bf 1.}] 
There exist two closed submodules $M_1$ and $M_2$ of $M$ such that $M_1$ is an essential extension of $X$, $M_1\cap M_2=0$ and $M$ is an essential extension of $M_1\oplus M_2$. In this case, $M_1$ and $M_2$ are, respectively, called a \textbf{closure} and  an \textbf{additive complement} of $X$ in $M$. 
\item[{\bf 2.}] 
If $M_3$ is any closed essential submodule of $M$ and $X\subseteq M_3$, then $M_3/X$ is an essential submodule of $M/X$.
\item[{\bf 3.}] 
If $M_4$ is a closed submodule of $M$ and $X\subseteq M_4$, then $M_4/M_1$ is a closed submodule of $M/X$.
\item[{\bf 4.}] 
If $M_1$ is a closed submodule of $M$ and $M_5$ is an essential submodule of $M$, then $M_1\cap M_5$ is a closed submodule of $M_5$.
\item[{\bf 5.}] 
If the kernel of a module homomorphism is essential, its image is singular . Consequently, if the quotient module $M/X$ is non-singular, then the submodule $X$ is closed in $M$.
\item[{\bf 6.}] 
If $M$ is non-singular, then the module $M/X$ is singular if and only if $M$ is an essential extension of $M_1$.
\item[{\bf 7.}] Let $A$ be a ring which has a minimal non-zero right ideal $S$. Then either $S^2=0$ or $S=eA$ for some idempotent $e$. In particular, if $A$ is semiprime ring, then $S=eA$ for some idempotent $e$.
\end{enumerate}

\textbf{Proof.} All assertions of Proposition 2.4  are well known. For example,  assertion \textbf{1} is verified via Zorn's lemma. Assertions \textbf{2}, \textbf{3}, \textbf{4} \textbf{5}, and \textbf{6} see  in  \cite[Proposition 1.4]{Goo76}, \cite[p.20, Ex.16]{Goo76} and  \cite[p.20, Ex.17]{Goo76}, \cite[Proposition 1.20]{Goo76}, \cite[Proposition 1.21]{Goo76}, respectively.  Assertion \textbf{1} is verified directly.\hfill$\square$

\textbf{Proposition 2.5.} Let $A$ be a right finite-dimensional semiprime ring. Assume that $A$ is either right non-singular or a right Goldie  ring (for a right finite-dimensional semiprime ring these assumptions are equivalent; e.g. see \cite[Corollary 3.32]{Goo76}).
\begin{enumerate}
\item[{\bf 1.}] 
The set of all essential right ideals of the ring $A$ coincides with the 
set of all right ideals of $A$ that contain regular elements. Consequently, the class of all singular (resp., non-singular) right $A$-modules coincides with the class of all torsion right $A$-modules (resp., torsion-free $A$-modules). All essential extensions of singular (resp., non-singular) right $A$-modules are singular (resp., non-singular) modules.
For any module $M_A$, the module $M/\text{Sing}(M)$ is non-singular. 
The ring $A$ has a semisimple artinian right classical ring of fractions $Q$, and non-zero injective non-singular uniform right $A$-modules coincide (up to an $A$-module isomorphism) with the minimal right ideals of the semisimple artinian ring $Q$; if the ring $A$ is prime, then  $Q$ is a simple ring isomorphic to a full matrix ring over a division ring. 
\item[{\bf 2.}]
Every non-torsion right $A$-module contains a non-zero non-singular submodule.
\item[{\bf 3.}]
Every non-singular divisible right $A$-module is injective.
\end{enumerate}

\textbf{Proof.} All assertions of Proposition 2.5  are well known;  see \cite{Goo76} and \cite{GooW89}. For example, assertion \textbf{1} is proved in \cite[5.9,~5.10,~6.14,~6.10(a)]{GooW89}, assertion \textbf{2} follows from \textbf{1}, assertion \textbf{3} is proved in \cite[6.12]{GooW89}.\hfill$\square$

\textbf{Proposition 2.6 \cite[Lemmas 4 and 5]{Tug26}.}\\ Let $A$ be a right Goldie prime ring  of cardinality $|A|$. 
\begin{enumerate}
\item[{\bf 1.}] 
Every non-zero ideal $B$ of $A$ is an essential right ideal and contains a regular element. 
\item[{\bf 2.}] 
There exists a positive integer $n$ such that for any non-zero elements $b_1,\dots\,b_n$ of the ring $A$, the module $b_1A\oplus\dots\oplus b_nA$ contains an isomorphic copy of the free cyclic module $A_A$.
\item[{\bf 3.}] 
If $X$ is a non-torsion right $A$-module and $X$ contains a non-singular submodule $Y$ of infinite Goldie dimension $\tau$, then $X$ contains a non-zero free submodule $F$ of infinite rank $\tau$; therefore, there exists an epimorphism from $F$ onto any $\tau$-generated right $A$-module.
\item[{\bf 4.}]  If $|A|$ is infinite and  $M$ is a non-singular right $A$-module of Goldie dimension $\tau\ge |A|$, then $M$ is not a Bassian module.
\end{enumerate} 

\textbf{Proposition 2.7.} Let $A$ be a ring such that there is a simple module $S_A$ which is not singular.
\begin{enumerate}
\item[{\bf 1.}] 
The simple module $S$ is nonsingular and is isomorphic to a nonsingular minimal right ideal $S'$ of $A$.
\item[{\bf 2.}] 
There are an ideal $B$ and a closed right ideal $C$ of $A$ such that $B\cap C=0$, $B\oplus C$ is an essential right ideal of $A$, $CB=0$, and $B_A$ is a non-singular semisimple non-zero module which is isomorphic to a direct sum of copies of the simple module $S$.
\item[{\bf 3.}] 
If $A$ is a right finite-dimensional prime ring, then $A$ is simple artinian.
\end{enumerate}

\textbf{Proof.} \textbf{1.} For the simple module $S_A$, there is an epimorphism $h\colon A_A\to S$ with kernel $M$, where $M$ is a maximal right ideal of $A$. Since the simple module $S$ is not singular, $A_A$ is not an essential extension of $M_A$. Therefore, $M\cap S'=0$ for some non-zero right ideal $S'$. Since $M$ is a maximal right ideal of $A$, we have $A_A=S'\oplus M$, $S'\cong A/M\cong S$, and $S'$ a non-singular minimal right ideal which is isomorphic to $S$.

\textbf{2.} The ring $A$ has a semisimple right ideal $B$ which is the sum of all minimal right ideals that are isomorphic to $S$. By Proposition 2.4(1), there is a closed right ideal  $C$ of $A$ such that $B\cap C=0$ and $B\oplus C$ is an essential right ideal of $A$. Since $B$ is the sum of all minimal right ideals that are isomorphic to the minimal right ideal $S'$, it is directly verified  that $B$ is an ideal of $A$.

\textbf{3.} Let $B$ and $C$ be the ideal and the right ideal from \textbf{2.}  By Proposition 2.6(1), the non-zero ideal $B$ is an essential right ideal and contains a regular element $b$. Since the module $B_A$ is semisimple, the cyclic right $A$-module $bA$  is semisimple. In addition, $A_A\cong bA$. Therefore,  $A$ is a prime semisimple ring, whence $A$ is a simple artinian ring.\hfill$\square$ 

\textbf{Proposition 2.8.} Let $A$ be a ring.
\begin{enumerate}
\item[{\bf 1.}] 
If all simple  $A$-modules are singular, then every generalized bassian non-singular right  $A$-module $M$ is bassian.
\item[{\bf 2.}]   
If $A$ is a right finite-dimensional prime non-artinian ring, then all generalized bassian non-singular right $A$-modules are bassian.
\end{enumerate}

\textbf{Proof.} \textbf{1.} Since  all simple  $A$-modules are singular, the non-singular module $M$ has no simple submodules. By Proposition 2.3(1),~ $M$ is bassian.

\textbf{2.} It follows from Proposition 2.7(3) that  all simple  $A$-modules are singular. By \textbf{1},  all generalized bassian non-singular right $A$-modules are bassian.\hfill$\square$ 

A module $A$ is \textbf{directly bounded} if there exists a finite $n$ such that for each simple submodule $T$ of $A$, there exists no direct summand of $A$ which is isomorphic to $T^{(n)}$.

\textbf{Proposition 2.9 \cite[Theorem 4.1]{DasK26}.} Let $M$ be a bassian module and let $S$ be a semisimple module. If the module $M$ is  directly bounded, then the module  $M\oplus S$ is generalized bassian. In particular, if $M$ is  finite-dimensional, then the module  $M\oplus S$ is generalized bassian.

\textbf{Proposition 2.10 \cite[Proposition 3.7]{DasK26}.} Let $A$ be a semiprime right Goldie ring, $Q$ be a classical right ring of fractions of $A$, and let $M$ be a right $A$-module such that  the (right) $Q$-module $M\otimes_AQ$ is finite-dimensional. Then $M\otimes_AQ$ is a bassian $Q$-module.

\textbf{Proposition 2.11 \cite[4.19, 4.20]{GooW89}.} \label{prop2}
A ring $A$ is right noetherian if and only if every injective right $A$-module is a direct sum of uniform injective modules, and if and only if all direct sums of injective right $A$-modules are injective.
\section{Modules over Non-Primitive HNP Rings}\label{sec3}
\textbf{Proposition 3.1 \cite[Proposition 4]{Tug26}.}
Let $A$ be a non-primitive HNP ring, $M$ be an injective indecomposable non-zero singular module, and let $M(P)$ be the primary component of $M$, where $P$ is a maximal invertible ideal of $A$. 
\begin{enumerate}
\item[{\bf 1.}] 
$M$ is a uniserial non-cyclic primary module without maximal submodules. All proper submodules of the module $M$ are cyclic modules of finite length. They form a countable chain $0=X_0\subset X_1\subset \ldots \subset X_k\subset \ldots$, where $X_k/X_{k-1}$ is a simple module for any $k$. There exists a positive integer $n$ such that $X_j/X_{j-1}\cong X_k/X_{k-1}$ if and only if $j-k$ is a multiple of $n$. In addition, the module $M$ has a surjective endomorphism with non-zero kernel. Therefore, $M$ is not Bassian.
\item[{\bf 2.}] 
If $\overline{M}$ is any non-zero homomorphic image of the module $M$, then for any cyclic submodule $\overline{X}$ of the module $\overline{M}$ of length $\ge k+n$, there exists an epimorphism $\overline{X}\to X_k$ with a non-zero kernel for an arbitrary cyclic submodule $X_k$ from \textbf{1}.
\item[{\bf 3.}] 
If $Y$ is an arbitrary cyclic uniserial $P$-primary module, then it is isomorphic to a subfactor of the module $M$ and is annihilated by some power of the ideal $P$.
\item[{\bf 4.}] 
If $\oplus_{i\in I}Y_i$ is a $P$-primary module, where all $Y_i$ are cyclic uniserial modules whose composition lengths have a common upper bound, then the module is annihilated by some positive power of the ideal $P$.
\item[{\bf 5.}] 
Every finite-dimensional reduced singular $A$-module is a finite direct sum of cyclic uniserial modules of finite composition length.
\item[{\bf 6.}] 
Every $A$-module with non-zero annihilator is a direct sum of cyclic uniserial modules of finite lengths with a common upper bound.
\end{enumerate}

\textbf{Proposition 3.2}  
Let $H$ be a generalized bassian module with a submodule $K\ne 0$.
\begin{enumerate}
\item[{\bf 1.}] 
If there is a monomorphism $\mu\colon H \to H/K$, then  there exists a semisimple submodule $S$ of $H$ and a direct sum decomposition  $H = S\oplus N$ such that $K$ is a semisimple  submodule of $S$ and $\mu(H)$ is an essential submodule of $N$.
\item[{\bf 2.}] 
Let $T$ be an infinite set and let  $K_t\cong K$ for all $t\in T$. If $H=\oplus_{t\in T}K_t$, then the  module  $K$ is semisimple.
\end{enumerate}
\textbf{Proof.} \textbf{1.} The assertion is proved in  \cite[Proposition 4.3]{DasK26}.

\textbf{2.}  We fix $t\in T$. Since the set $T$ is infinite, there exists an epimorphism $H\to \oplus_{u\in T\setminus t}K_t$ with kernel $K_t$. Therefore,  there is a monomorphism $\mu\colon H \to H/K$. By \textbf{1}, the module $K$ is semisimple.\hfill$\square$

\textbf{Proposition 3.3.}
Let $A$ be a non-primitive HNP ring, $M$ be a singular right $A$-module, and let $\{M(P_i)\}$ be the set of all primary components of $M$. 
\begin{enumerate}
\item[{\bf 1.}] 
All finitely generated submodules of the module $M$ are finite direct sums of cyclic uniserial modules of finite length. 
\item[{\bf 2.}] 
Every non-injective submodule of the module $M$ has a non-zero cyclic uniserial direct summand of finite length.
\item[{\bf 3.}] 
$M$ is an injective module if and only if every primary component $M(P_i)$ is an injective module if and only if every primary component $M(P_i)$ is a direct sum of uniserial injective modules.
\item[{\bf 4.}] 
If $M$ is a generalized bassian module, then $M$ is a reduced module.
\item[{\bf 5.}] 
If $M$ is a generalized bassian module and all its homogeneous components are finite-dimensional, then $M$ is  a reduced Bassian module.
\item[{\bf 6.}] 
The module $M$ is generalized bassian if and only if all primary components of $M$ are  generalized bassian modules.
\item[{\bf 7.}] 
If every primary component of $M$ is a direct sum of a Noetherian module and a semisimple module, then $M$ is generalized bassian.
\item[{\bf 8.}] 
A submodule $B$ of $M$ is called a \textbf{basic submodule} $B$ if $B$ is a direct sum of cyclic uniserial modules of finite composition length,  $B$ is a pure submodule of $M$, the quotient module $M/B$ is a direct sum of uniserial non-cyclic modules $M_i/B$, and  any two basic submodules of the module $M$ are isomorphic; it follows from Propositions 3.1(1) and  2.11 that each of the modules $M_i/B$ and $M/B$ is injective uniserial.\\ 
A basic submodule $B$ of $M$ is called a \textbf{lower basic submodule} if $\text{u.dim }(M/B)=\text{u.dim }M$. If  the set $I$ is finite, then every basic submodule of $M$ contains  a lower basic submodule.\\
There is a direct sum decomposition $M=H_1\oplus H_2$, where $H_1$ has a non-zero annihilator $r(H_1)$ and $\text{u.dim }H_2 = \text{u.dim}(H_2/B_2) = \text{u.dim}(M/B)$, where $B$ and $B_2$ are lowest basic submodules of $M$ and $H_2$, respectively.
\item[{\bf 9.}] 
If  $M$ is a generalized bassian primary module, then $M$ has a non-zero annihilator $r(M)$,  and there is a positive integer $n$ such that  $M=\oplus_{j\in J}X_j$, where every $X_j$ is a uniserial noetherian cyclic module of length $\le n$.
\end{enumerate}

\textbf{Proof.} \textbf{1, 2, 3.} The assertions are well known; e.g., see  \cite[Proposition 3]{Tug26}. Also, see  \cite{Sin75} and \cite{Sin76}.

\textbf{4.} Assume that  $M$ is not a reduced module. Then $M$  has a non-zero injective direct summand. By \textbf{3}, $M$  has a non-zero injective uniserial direct summand $M'$. By Proposition 2.2, the module $M'$ is generalized Bassian. By Proposition 2.3(1), the  module $M'$ is  Bassian. This contradicts to Proposition 3.1(1).

\textbf{5.}  By \textbf{4}, the module $M$ is reduced. Since  all homogeneous components of $M$ are finite-dimensional,  every homogeneous component of $\text{Soc }(M)$ is finitely generated. By Proposition 2.3(1), the module $M$ is Bassian.

\textbf{6.} If  $M$ is generalized bassian, then all primary components of $M$ are  generalized bassian by Proposition 2.2. 

Let  all primary components $M_i$ of $M$ be generalized bassian modules. We have $M=\oplus_{i\in I}M_i$. By the use of Remark 1.4, we may directly verify that $M$ is a  generalized bassian module.

\textbf{7.} The assertion follows from Proposition 2.9 and Remark 2.1.

\textbf{8.} For the first assertion, see  \cite[Theorem 1]{Sin78}; see also \cite[Theorem 1.3]{BenS83}. For the second and third assertions, see  \cite[Theorems 3.8, 3.10]{SinT79}.

\textbf{9.} By \textbf{8}, there is a direct decomposition $M=H_1\oplus H_2$, where $H_1$ has a non-zero annihilator $r(H_1)$ and $\text{u.dim }H_2 = \text{u.dim}(H_2/B_2)$, where $B_2$ is a lowest basic submodule of $H_2$. Since $r(H_1)\ne 0$, the quotient ring $A/r(H_1)$ is a serial Artinian ring. Therefore,  there is a positive integer $n$ such that  $M=\oplus_{j\in J}X_j$, where every $X_j$ is a uniserial noetherian cyclic module of length $\le n$. By Proposition 2.2, the module $H_2$ is generalized bassian. 

It remains to show that $H_2=0$. Assume that $H_2\ne 0$. Then $B_2\ne 0$,~ $H_2$ has the socle $S_2\ne 0$, and $H_2$ is an essential extension of $S_2$. Since  $\text{u.dim }H_2 = \text{u.dim}(H_2/B_2)$ and the module $H_2/B_2$ is injective, there exists a monomorphism $f\colon S\to H_2/B_2$, which can be extended to a homomorphism $\mu$ from $H_2$ into the injective module $H_2/B_2$. Since $S\cap\text{Ker }g=\text{Ker }f=0$ and $H_2$ is an essential extension of $S_2$, we have that $\mu$ is a monomorphism from the generalized bassian module $H_2$ to $H_2/B_2$, where $B_2\ne 0.$ 

Now we use Proposition 3.2, where we set  $H=H_2$ and $K=B_2$. Since  $\mu\colon H\to H/K$ is a monomorphism, it follows from Proposition 3.2 that  there exists a semisimple submodule $S$ of $H$ and a direct decomposition  $H = S\oplus N$ such that $K$ is a semisimple  submodule of $S$ and $\mu(H)$ is an essential submodule of $N\ne 0$. Since $S$ is a semisimple module, there are direct decompositions $S=K\oplus S'$ and $H = K\oplus S'\oplus N$, Therefore, the non-zero module $N$ is isomorphic to a direct summand of the non-zero injective module $H/K$. In addition, $N$ is  generalized bassian, since $H$ is  generalized bassian.  Thus, there exists a   generalized bassian non-zero injective singular module $N$. By Propositions 2.11 and 2.2, there exists a  generalized bassian non-zero injective indecomposable singular module $N'$ which is  bassian, by Proposition 2.3(2). This contradicts to Proposition 3.1(1).\hfill$\square$

\textbf{Proposition 3.4 \cite[Proposition 5]{Tug26}.}
Let $A$ be a non-primitive Dedekind prime ring, $M_A$ be an injective indecomposable non-zero singular module, and let $M(P)$ be the primary component of $M$, where $P$ is a maximal ideal. 
\begin{enumerate}
\item[{\bf 1.}] 
$M$ is a uniserial non-cyclic primary module, all non-zero quotient modules of $M$ are isomorphic to $M$, all proper submodules of $M$ are cyclic modules of finite length and form a countable chain $0=X_0\subset X_1\subset \ldots \subset X_k\subset \ldots$, where for any $k\ge 1$, all modules $X_k/X_{k-1}$ are isomorphic simple modules. 
\item[{\bf 2.}] 
If $Y$ is an arbitrary cyclic uniserial $P$-primary module, then it is isomorphic to a subfactor of the module $M$ and is annihilated by some power of the ideal $P$.
\item[{\bf 3.}] 
Every $P$-primary module $N_A$ of Goldie dimension $\tau$ is an essential extension of a direct sum $\oplus_{j\in J}S_j$ of isomorphic simple modules $S_j$, where $|J|=\tau$ and $(\oplus_{j\in J}S_j)P=0$.
\item[{\bf 4.}] 
If $\oplus_{i\in I}Y_i$ is a $P$-primary module, where all $Y_i$ are cyclic uniserial modules whose composition lengths have a common upper bound, then the module is annihilated by some positive power of the ideal $P$ and  all simple subfactors of $\oplus_{i\in I}Y_i$ are isomorphic.
\item[{\bf 5.}] 
Let $Q$ be an injective $P$-primary $A$-module of Goldie dimension $\tau$ and let $X$ be a $P$-primary $A$-module of Goldie dimension $\sigma\le\tau$. Then there exists a monomorphism $f\colon X\to Q$.
\end{enumerate}

\textbf{Proposition 3.5 \cite[Theorem 1]{Tug26}. } If $A$ is a non-primitive Dedekind prime ring and $M$ is a singular right $A$-module with primary components $M_i$, $i\in I$, then the following conditions are equivalent. 
\begin{enumerate} 
\item[ 1)]  
$M$ is a Bassian module. 
\item[2)] 
Every primary component $M_i$ of the module $M$ is a reduced finite-dimensional module. \item[3)] 
Every primary component $M_i$ of the module $M$ is a noetherian module.
\end{enumerate}

\textbf{Proposition 3.6.} If $A$ is a non-primitive Dedekind prime ring and $M$ is a generalized bassian singular right $A$-module, then every primary component $M_i$ of $M$ is a direct sum of a noetherian module and a semisimple module.

\textbf{Proof.}  By Proposition 2.2, every primary component $M_i$ of $M$ is a  generalized bassian primary module. Without loss of generality, we assume that $M=M_i$ for some $i\in I$. By Proposition 3.3(9), there is a positive integer $n$ such that  $M=\oplus_{j\in J}K_j$, where every $K_j$ is a uniserial noetherian cyclic module of length $\le n$. By Proposition 3.4(4), all simple subfactors of $M$ are isomorphic. Let $I$ be the subset of $J$ such that the length of $K_i$  exceeds $2$. We denote the module $\oplus_{i\in I}K_i$ by $M_1$ and we set $M_1=0$ if the set $I$ is empty. We denote the semisimple module $\oplus_{j\in J\setminus I}K_j$ by $M_2$; if $I=J$, then we set $M_2=0$. We have $M=M_1\oplus M_2$. It remains to prove that the module $M_1$ is noetherian. This is true provided the set $I$ is finite or empty.

Assume that the set $I$ is infinite. By Proposition 2.2, the module $M_1$ is generalized bassian. Since all simple subfactors of every cyclic uniserial summand $K_i$ of length $2\le d\le n$ are isomorphic to a fixed simple module, there is an infinite subset $T$ of $I$ such that all the modules $X_t$, $t\in T$ are isomorphic. We set $H=\oplus_{t\in T}X_t$. By Proposition 2.2, the module $H$ is generalized bassian. By Proposition 3.2(2), the modules  $X_t$ are semisimple. This contradicts to the $X_t$ being of length $\ge 2$. \hfill$\square$
\section{Remarks}\label{sec4}
\textbf{Remark 4.1.} Theorem 1.2 is the main result of  the given paper and follows from Propositions 3.3(7) and 3.6.

A submodule $N$ of a module $M$ is said to be \textbf{fully invariant} if $f(N)\subseteq N$ for every endomorphism $f$ of $M$. A module is called an \textbf{invariant} module provided every its submodule is fully invariant. A ring $A$ is right invariant  if and only if every right ideal of $A$ is an ideal of $A$; clearly, commutative rings and division rings are right and left invariant rings.
 
\textbf{Remark 4.2 \cite[Examples 2.1 and 2.7]{DasK26}.}\\ Any prime or right invariant ring $A$ is a bassian right $A$-module.

A submodule $N$ of a module $M$ is said to be \textbf{dense} if, for any $y\in M$ and non-zero $x\in M$, there exists $r\in R$ such that $xr = 0$, and $yr\in N$. 
A module $M$ is said to be \textbf{monoform} if every nonzero homomorphism from a submodule $N$ of $M$ to $M$ is a monomorphism. A module $M$ is said to be \textbf{polyform} if $\text{Ker}(f)$ is a closed submodule of $K$ for every submodule $K$ of $M$ and $f\colon K \to M$.

\textbf{Remark 4.3 \cite[Examples 2.11 and 2.5, Propositions 3.2 and 3.3]{DasK26}.}\\ 
Let $M$ be a module.
\begin{enumerate}
\item[{\bf 1.}] 
If all proper submodules of $M$ are dense, then $M$ is bassian.
\item[{\bf 2.}] 
If  $M$ is monoform, then $M$ is bassian.
\item[{\bf 3.}] 
If $M$ is uniform and $M$ is nonsingular or polyform, then $M$ is bassian.
\end{enumerate}

\textbf{Remark 4.4 \cite[Corollary 3.2]{DasK26}.} Let $B$ be an ideal of a ring $A$ and $M$ be a right $A$-module. If $M/MB$ is a bassian $R/B$-module, then $M$ is a bassian right $A$-module.

\textbf{Remark 4.5 \cite[Proposition 3.8]{DasK26}.} If $M$ is an $A$-module such that $M[X]$ is a bassian $A[X]$-module, then the $R$-module $M$ is bassian.

\textbf{Remark 4.6. Open question.} Can Theorem 1.2 be extended from Dedekind prime rings to HNP rings?

Is Theorem 1.2 true for HNP rings?

\end{document}